\documentclass[10pt]{article}
\usepackage{amsmath,amsthm,amsfonts}
\usepackage{graphicx}
\usepackage{hyperref}
\usepackage[usenames]{color}

\newtheorem{theorem}{Theorem}[section]
\newtheorem{lemma}[theorem]{Lemma}
\theoremstyle{definition}

\newtheorem{corollary}[theorem]{Corollary}

\newcounter{comcount}

\title{Twisted conjugacy classes in nilpotent groups}

\author{V. ROMAN'KOV\footnote{The  author was partially supported by RFBR, Grant 07-01-00392}
}

\begin{document}

\maketitle

\begin{abstract}
Let $N$ be a finitely generated nilpotent group. Algorithm is
constructed such, that for every automorphism $\varphi \in Aut(N)$
defines the Reidmeister number $R(\varphi ).$ It is proved that
any free nilpotent group of rank $r = 2$ or $3$  and class $c \geq
4r,$ or rank $r \geq 4$ and class $c \geq 2r,$ belongs to the
class $R_{\infty }.$
\end{abstract}

\tableofcontents

\medskip
\begin{quote}{\it \scriptsize 2000 Mathematics Subject Classification.
Primary 20F10. Secondary 20F18; 20E45; 20E36. }

\medskip
\noindent {\it \scriptsize Keywords.} {\rm \scriptsize nilpotent
groups, twisted conjugacy classes, Reidemeister numbers,
automorphisms.}
\end{quote}

\bigskip

\section{Introduction}
\label{se:intro}

Let $\varphi : G \rightarrow G$ be an automorphism of a group $G.$
One says that the elements $g, f \in G$ are $\varphi -${\it
twisted conjugated}, denoted by $g \sim_{\varphi } f,$ if and only
if there exists $x \in G$ such that $g = (x\varphi )^{-1}fx,$ or
equivalently $(x\varphi )g = fx.$ A class of equivalence
$[g]_{\varphi }$ is called the {\it Reidemeister class} (or the
$\varphi -${\it conjugacy class of } $\varphi $). The number
$R(\varphi )$ of Reidemeister classes is called the {\it
Reidemeister number of} $\varphi .$

There are different origins of interest in twisted conjugacy
relations. In classical Nielsen-Reidemeister fixed point theory,
$R(\varphi )$ plays a crucial role in estimating the Nielsen
number $N(f)$ of a selfmap $f : C \rightarrow C,$ where $C$ is a
compact connected manifold. In such setting $\varphi $ appears as
the homomorphism induced by $f$ on the fundamental group
$\pi_1(C).$ The Nielsen number $N(f)$ is a homotopy invariant
which provides a lower bound of size of the set $Fix_f(C)$ of
fixed $f-$points. The Selberg theory (see \cite{Sel}), and
Algebraic Geometry (see \cite{AG}) present other sources for the
twisted conjugacy relations.

One of the central problems in the field concerns obtaining a
twisted analogue of the classical Burnside-Frobenius theorem, that
to show the coincidence of the Reidemeister number $R(\varphi )$
and the number of fixed points of the induced homeomorphism of an
appropriate dual object. The authors of the paper \cite{FLT}
emphasize that one step in this process is to describe the class
of groups $G,$ such that $R(\varphi ) = \infty $ for any
automorphism $\varphi : G \rightarrow G.$ In a number of papers
(\cite{TW}, \cite{FGW}, \cite{FLT}) this class of groups is
denoted $R_{\infty }.$ Namely, a group $G$ has {\it property}
$R_{\infty }$ (or is an $R_{\infty }$ {\it group}) if all of its
automorphisms $\varphi $ have $R(\varphi ) = \infty .$

It was shown by various authors that the following groups belong
to the class $R_{\infty }:$

\begin{itemize}
\item non-elementary Gromov hyperbolic groups \cite{LL, Fe}; \item
Baumslag-Soliter groups $B(m, n)$ except for $B(1, 1)$ \cite{FG2};
\item generalized Baumslag-Soliter groups, that is finitely
generated groups which act on a tree with all edge and vertex
stabilizers infinite cyclic \cite{L, TWo2}; \item lamplighter
groups ${\bf Z}_n wr {\bf Z}$ if and only if $2|n$ or $3|n$
\cite{GW1}; \item the solvable generalization $\Gamma $ of $BS(1,
n)$ given by the short sequence $1 \rightarrow {\bf Z}[1/n]
\rightarrow \Gamma \rightarrow {\bf Z}^k \rightarrow 1 $
\cite{TWo}; \item relatively hyperbolic groups \cite{FG1}; \item
the Grigorchuk group and Gupta-Sidki group \cite{FLT}.
\end{itemize}

For the immediate consequences of $R_{\infty }$ property for
topological fixed point theory see, e.g. \cite{TWo2}.

In the present paper we treat with finitely generated  nilpotent
groups. On the first glance  nilpotent groups are too far to
possess the $R_{\infty }$ property. Indeed, E.G. Kukina \cite{Kuk}
noted that the Reidemeister spectrum

$$Spec_R(G) = \{ R(\varphi ) | \varphi \in Aut(G)\}$$

\noindent of any free abelian group $G = {\bf Z}^k, \  k \geq 2,$
coincides with ${\bf N} \cup \{\infty \},$ i.e. is full.

She also calculated the Reidemeister spectrums of the free
nilpotent groups $N_{rc}$ of rank $r$ and class $c$ in cases $r =
2, 3$ and $c = 2.$ Namely,

$$Spec_R(N_{22}) = \{ 2k | k \in {\bf N}\} \cup \{ \infty \},$$

$$Spec_R(N_{32}) = \{ 2k+1 | k \in {\bf N} \} \cup \{ 4k | k \in
{\bf N}\} \cup \{\infty \}.$$

Note that the occurrence of any even number in  $Spec_R(N_{22})$
was established by F.K. Indukaev \cite{Ind}.

\bigskip

\section{The representatives of the twisted conjugacy classes}
\label{se:rep}

Let $G$ be any finitely generated group, and let $C$ be a central
subgroup of $G.$ For any automorphism $\varphi : G \rightarrow G$
define a central subgroup

\begin{equation}
\label{eq:1}L(C, \varphi ) = \{ c \in C | \exists x \in G |
x\varphi = cx\}.
\end{equation}

\begin{lemma}
\label{le:1} Any pair of elements $c_1, c_2 \in C$ are $\varphi
-$conjugated in $G$ if and only if
\begin{equation}
\label{eq:2} c_1^{-1}c_2 \in L(C, \varphi ).
\end{equation}
\end{lemma}

Proof. Let $c_1 \sim_{\varphi } c_2 \Rightarrow \exists x \in G |
(x\varphi ) c_1 = c_2 x \Rightarrow x\varphi = c_1^{-1}c_2 x
\Rightarrow c_1^{-1}c_2 \in L(C, \varphi ).$

In opposite, let $c_1^{-1}c_2 \in L(C, \varphi ) \Rightarrow
\exists x \in G | x\varphi = c_1^{-1}c_2 x \Rightarrow (x\varphi
)c_1 = c_2 x \Rightarrow c_1 \sim_{\varphi } c_2.$

\begin{corollary}
\label{co:2} The set of all $\varphi -$conjugacy classes of the
elements $c \in C$ coincides with the set of all cosets $C$ w.r.t.
$L(C, \varphi ).$ In particular, $R(\varphi ) \geq [C : L(C,
\varphi )].$
\end{corollary}

\begin{corollary}
\label{co:3} In the case of a finitely generated abelian group $G
= C$ there is an effective procedure for calculating the
Reidemeister number $R(\varphi )$ of any automorphism $\varphi : G
\rightarrow G.$ Moreover, such procedure finds in a finite case a
set of representatives of all $\varphi -$conjugacy classes.
\end{corollary}

Proof. Let $G = gp(a_1, ..., a_p), $ and $a_i\varphi = c_ia_i, \ i
= 1, ..., p.$ Then $L(C, \varphi ) = gp(c_1, ..., c_p).$ A
standard procedure easily defines either $[C : L(C, \varphi )] =
\infty ,$ or not.

In a finite case we just need to calculate the index $[C : L(C,
\varphi )]$ and to find any transversal set of $C$ w.r.t. $L(C,
\varphi ).$ It is well known that such effective procedure exists
in any finitely generated effective presented abelian group $G.$

Now let $G$ be any group, and $\varphi : G \rightarrow G$ be an
automorphism of $G.$ Let $C$ be a central $\varphi -$admissible
subgroup of $G.$ Denote by $\eta : G \rightarrow \bar{G} = G/C$
the standard epimorphism, and by $\bar{\varphi } : \bar{G}
\rightarrow \bar{G}$ the induced by $\varphi $ automorphism. Let
$[\bar{g}]_{\bar{\varphi }}$ be any $\bar{\varphi }-$conjugacy
class in $\bar{G},$ and $g \in G$ be any pre image of $\bar{g}$ in
$G.$ Define an automorphism $\varphi _g : G \rightarrow G$ by
$\varphi _g = \varphi \circ \sigma _g,$ where $\sigma_g \in InnG,
\sigma_g : h \mapsto g^{-1}hg$ for all $h \in G.$ Define as before
a subgroup $L(C, \varphi_g)$ of $G.$ Then we have

\begin{lemma}
\label{le:4} The full pre image in $G$ under the standard
epimorphism $\eta : G \rightarrow \bar{G}$ of any $\bar{\varphi
}-$equivalence class $[\bar{g}]_{\bar{\varphi }}$ is a disjoint
union of $s = [C : L(C, \varphi_g)]$ $\varphi -$conjugacy classes
\begin{equation}
\label{eq:3} [gc_1]_{\varphi} \cup ... \cup [gc_s]_{\varphi },
\end{equation}

\noindent where $\{c_1, ..., c_s\}$ is any transversal set of $C$
w.r.t. $L(C, \varphi_g).$
\end{lemma}

Proof. Let we show first that all $\varphi -$conjugacy classes in
(\ref{eq:3}) are different. Suppose that $gc_i \sim_{\varphi }
gc_j$ for some $i \not= j.$ Then

\begin{equation}
\label{eq:4} \exists x \in G | (x\varphi )gc_i = gc_jx \Rightarrow
(x\varphi_g) = c_i^{-1}c_jx \Rightarrow c_i^{-1}c_j \in L(C,
\varphi_g),
\end{equation}

\noindent which contradicts to our assumption.

Now we are to prove that every element $gc, c \in C,$ is $\varphi
-$conjugated to some of $gc_i, i = 1, ..., c_s.$ Namely, if $c =
c_ic',$ where $c' \in L(C, \varphi_g )$ then $gc \sim_{\varphi }
gc_i.$ To prove it suppose that for $x \in G$ one has $x\varphi_g
= c'x.$ Then
\begin{equation}
\label{eq:5} x\varphi = g(x\varphi_g)g^{-1} = gxg^{-1}c'
\Rightarrow (x\varphi )gc_i = gc_ic'x \Rightarrow gc \sim_{\varphi
} gc_i.
\end{equation}.

At last suppose that any element $f \in G$ belongs to the full pre
image of $[\bar{g}]_{\bar{\varphi }}.$ So,

\begin{equation}
\label{eq:6} \bar{f} \sim_{\bar{\varphi }}\bar{g} \Rightarrow
\exists x \in G, c \in C | (x\varphi )g = fxc \Rightarrow
(x\varphi )gc^{-1} = fx \Rightarrow f \sim_{\varphi } gc_i,
\end{equation}

\noindent where $c^{-1} = c_ic', c' \in L(C, \varphi_g).$

The following result was obtained in \cite{RV} even in more
general form.

\begin{lemma}
\label{le:5} Let $N$ be a finitely generated nilpotent group of
class $k,$ and $\varphi : N \rightarrow N$ be any automorphism.
Then there is an effective procedure which gives a finite
generating set of a subgroup
\begin{equation}
\label{eq:7} Fix_{\varphi }(N) = \{ x \in N | x\varphi = x \}.
\end{equation}
\end{lemma}

Proof. Obviously, such procedure exists in abelian case $k = 1.$

Let $C$ be a central $\varphi -$admissible subgroup of $N,$ one
can take the last nontrivial member $C = \gamma_{k-1}N$ of the
lower central series  of $N.$ Let $\bar{\varphi } : N/C
\rightarrow N/C$ be the induced by $\varphi $ automorphism. By
induction on $k$ we have a generating set of
\begin{equation}
\label{eq:8} Fix_{\bar{\varphi }}(N/C) = gp(\bar{g_1}, ...,
\bar{g_l}).
\end{equation}

Then we can effectively obtain a generating set $\{g_1, ..., g_l,
g_{l+1}, ..., g_p\}$ of the full pre image of $Fix_{\bar{\varphi
}}(N/C)$ in $N.$ Namely, we take some pre images $g_1, ..., g_l$
of the elements $\bar{g_1}, ..., \bar{g_l},$ respectively, and add
a generating set $\{g_{l+1}, ..., g_p\}$ of $C.$ See for details
\cite{KRRRC}. So,
\begin{equation}
\label{eq:9}
 g_i\varphi = c_ig_i, \  i = 1, ..., p.
\end{equation}

Introduce a homomorphism of the full pre image of
$Fix_{\bar{\varphi }}(N/C)$ in $N$ to $C,$ uniquely defined by the
map
\begin{equation}
\label{eq:10} \mu : g_i  \mapsto  c_i, i = 1, ..., p.
\end{equation}

Since the derived subgroup belongs to $ker\mu $ we actually have a
homomorphism of abelian groups. Then we can effectively find a
generating set of
\begin{equation}
\label{eq:11} Fix_{\varphi }(N) = ker\mu .
\end{equation}

\bigskip
{\bf Theorem 1.} Let $N$ be a finitely generated nilpotent group
of class $k.$ Then there is an effective procedure to calculate
for any automorphism $\varphi : N \rightarrow N$ the Reidemeister
number $R(\varphi ).$

Proof. By induction on $k,$ starting with Corollary \ref{co:3}, we
can assume that the statement is true for any finitely generated
nilpotent group of class $\leq k-1.$

Let $\varphi : N \rightarrow N$ be any automorphism of $N.$
Consider some $\varphi -$admissible central series
\begin{equation}
\label{eq:12} N = C_1 > C_{2} > ... > C_k > 1.
\end{equation}

For example, one can take the lower central series in $N.$ The
automorphism $\varphi $ induces the automorphisms $\varphi_i : N_i
\rightarrow N_i, $ where $N_i = N/C_{i+1}, \  i = 1, ..., k-1.$

By induction we have the number $R(\varphi _{k-1}).$ Obviously,
$R(\varphi_{k-1}) = \infty $ implies $R(\varphi_k) = R(\varphi ) =
\infty .$

Suppose, $R(\varphi_{k-1}) = r < \infty .$ Let
\begin{equation}
\label{eq:13} [\bar{g_1}]_{\varphi_{k-1}}, ...,
[\bar{g_{r}}]_{\varphi_{k-1}}
\end{equation}

\noindent be the set of all $\varphi_{k-1}-$conjugacy classes in
$N_{k-1} = N/C_k.$ As before $\bar{f}$ means the image of any
element $f \in N$ under the standard epimorphism $N \rightarrow
N_{k-1}.$

To apply Lemma \ref{le:4} we need in subgroups $L(C_k,
\varphi_{g_i}), \ i = 1, ..., r.$ Let $\psi : N \rightarrow N$ be
any automorphism with $\psi (C_k) = C_k,$ and $\psi_{k-1} :
N_{k-1} \rightarrow N_{k-1}$ be the automorphism induced by $\psi
.$

Firstly we derive by Lemma \ref{le:5} a generating set $\{f_1,
..., f_q\}$ of the full pre image  $H_{\psi }$ of
$Fix_{\psi_{k-1}}(N_{k-1})$ in $N.$ Then we compute

\begin{equation}
\label{eq:14} f_j\psi = c_jf_j, \  c_j \in C_k, \ j = 1, ..., q,
\end{equation}

\noindent and conclude that

\begin{equation}
\label{eq:15} L(C_k, \psi ) = gp(c_1, ..., c_q).
\end{equation}

We repeat such process for every $\psi = \varphi_{g_i}, \  i = 1,
..., r.$ Every time we apply Lemma \ref{le:4} to find a
presentation of the full pre image $H_{\varphi_{g_i}}$ of
$[\bar{g_i}]_{\varphi_{k-1}}$ in $N$ as a disjoint union of
$\varphi -$conjugacy classes in $N.$ The union of them is the set
of all $\varphi -$conjugacy classes in $N.$

Theorem is proved.

\bigskip

\section{Most of the free nilpotent groups are in $R_{\infty }$}
\label{se:rep}

Let $N$ be a finitely generated torsion free nilpotent group.
Every quotient $A_i = \zeta_{i+1}N/\zeta_{i}N, \ i = 0, 1, ..., k
- 1,$ of the upper central series

\begin{equation}
\label{eq:16} \zeta_oN = 1 < \zeta_1N < ... < \zeta_kN = N
\end{equation}

\noindent is a free abelian group of a finite rank. Any
automorphism $\varphi : N \rightarrow N$ induces the automorphisms
$\varphi_i : N_i \rightarrow N_i, $ where $N_i = N/\zeta_iN, \ i =
1, ..., k-1,$ as well as the automorphisms $\bar{\varphi_i} : A_i
\rightarrow A_i,$  \  i = 0, 1, ..., k - 1.

\begin{lemma}
\label{le:3.1} Let $A$ be a free abelian group of a finite rank
$r.$ Let $\varphi : A \rightarrow A$ be any automorphism of $A$
such that

\begin{equation}
\label{eq:17} Fix_{\varphi }A \not= 1.
\end{equation}

Then $R(\varphi ) = \infty .$ \end{lemma}

Proof. It was shown in Lemma \ref{le:1} that $R(\varphi ) = [A :
L(A, \varphi )],$ where in our case

\begin{equation}
\label{eq:18} L(A, \varphi ) = im (\varphi - id ).
\end{equation}

The subgroup $L(A, \varphi )$ is a free abelian group of rank $s =
range (\varphi - id ).$ If $Fix_{\varphi }(A) \not= 1$ we have $s
< r$ and so $[A : L(A, \varphi )] = \infty .$

\begin{corollary}
\label{co:3.2} Let $N$ be a finitely generated torsion free
nilpotent group of class $k,$ and $\varphi : N \rightarrow N$ be
any automorphism of $N.$ Suppose that $Fix_{\bar{\varphi_i}}(A_i)
\not= 1$ for some $i = 0, 1, ..., k-1.$ Then $R(\varphi ) = \infty
.$
\end{corollary}

Proof. We can assume that $i$ is maximal with property
$Fix_{\bar{\varphi_i}}(A_i) \not= 1.$ Then $Fix_{\varphi_{i+1}} =
1.$ Let in the denotions of Lemma \ref{le:1} $G = N_i$ and $C =
A_i.$ Then by the statement of Lemma \ref{le:3.1} $R(\bar{\varphi
_i}) = \infty ,$ so $R(\varphi _i) = \infty $ too.

E. Formanek \cite{F} classified all pairs $r, c$ for which the
free nilpotent group $N_{rc}$ of rank $r \geq 2$ and class $c \geq
2 $ has nontrivial elements fixed by all automorphisms. Note that
V. Bludov (see \cite{B}, background to Problem N1 by A. Myasnikov)
gave the first examples of such elements for $r = 2, c = 4k, k
\geq 2.$

\bigskip
{\bf Theorem of Formanek.} Let $N_{rc}$ be a free nilpotent group
of rank $r$ and class $c.$ Then there are nontrivial elements of
$N_{rc}$ which are fixed by all automorphisms of $N_{rc}$ if and
only if

(a) $r = 2$ or $r = 3,$  and $c = 2kr, \  k \geq 2.$

(b) $r \geq 4$ and $c = 2kr, \  k \geq 1.$

\bigskip
This remarkable result allows to define place of most free
nilpotent groups of finite rank in the class $R_{\infty }.$

\bigskip
{\bf Theorem 2.} Let $N_{rc}$ be a free nilpotent group of rank $r
\geq 2$ and class $c \geq 2.$ Suppose that for $r = 2$ or $r = 3$
we have $c \geq 4r$ and for $r \geq 4$ we have $c \geq 2r.$ Then
$N_{rc} \in R_{\infty }.$

Proof. Denote $N = N_{rc}.$ Let $\varphi : N \rightarrow N$ be any
automorphism. It follows from the Formanek's theorem that there is
a maximal $k$ such that the induced automorphism $\varphi_j : N_j
\rightarrow N_j$ has $Fix_{\varphi_j}(N_j) \not= 1.$ Then
$Fix_{\bar{\varphi_j}}(A_j) \not= 1,$ and so by Lemma \ref{le:1}
$R(\varphi_j) = \infty .$ By Corollary \ref{co:3.2} we have
$R(\varphi ) = \infty ,$ which gives the conclusion of the
theorem.

\end{document}